Dmytro Taranovsky
April 20, 2005

# Extending the Language of Set Theory

*Note (March 10, 2012)*: A follow-up paper on higher order theory is now available: "Reflective Cardinals", arXiv:1203.2270.

**Abstract:** We discuss the problems of incompleteness and inexpressibility. We introduce almost self-referential formulas, use them to extend set theory, and relate their expressive power to that of infinitary logic. We discuss the nature of proper classes. Finally, we introduce and axiomatize a powerful extension to set theory.

**Contents:** Introduction and Outline, Incompleteness and Inexpressibility, Almost Self-Referential Formulas, Expressive Power of the Extensions, A Hierarchy of Large Cardinals, Proper Classes, Higher Order Set Theory.

*Update (December 19, 2016)*: Added references and revised two sections ("Almost Self-Referential Formulas" and "Expressive Power of the Extensions"); content updates in other sections are marked with 'update'.

## Introduction and Outline

Despite vast advances in set theory and mathematics in general, the language of set theory, which is also the language of mathematics, has remained the same since the beginning of modern set theory and first order logic. That formal language has served us well, but it is of necessity limited, and does not adequately deal with properties for which there is no set of all objects satisfying the property. The purpose of this paper is to address the deficiencies in the expressive power of ordinary set theory.

We assume general familiarity with mathematical logic and set theory. A good exposition of set theory is in [1] or [2], the later having more emphasis on large cardinals, including history and philosophy.

The second section discusses incompleteness and inexpressibility in general. Much of the paper discusses logic that allows some self-reference but with limitations to prevent loops and infinite regress. The logic is more expressive than (a formulation of) infinitary logic. We show that levels of infinitary logic correspond to levels of almost self-referential logic and levels of constructible hierarchy above V. We use strong logics and weak set theories to clarify the strengths of inaccessible and Mahlo cardinals. We also use large cardinals to get reflection principles, and use reflection principles to axiomatize extensions to the language of set theory.

One section discusses several approaches to proper classes, and shows how these approaches benefit from the material in this paper. The last section (which is largely independent of previous sections) introduces the notion of reflective ordinals [update (Dec 16, 2016): which are now called reflective cardinals], which captures higher order set theory and is likely to supersede extensions to set theory proposed earlier in the paper. However, its theory is less certain, and constructs such as almost self-referential logic are important in their own right.

# Incompleteness and Inexpressibility

Introduction of formal languages and axiomatic systems into mathematics has created the necessary precision in the language of mathematics and its methods. However, every formal language fails to express some mathematical claims and every (recursive) sufficiently expressive axiomatic system is incomplete.

To avoid incompleteness, some had proposed to strengthen the deductive apparatus. However, if the new deductive apparatus is not recursive, then to be useable by mathematicians, a proof in the new deductive system would have to include a demonstration that it is a valid proof. To achieve the necessary precision and to avoid basic disagreements that plague many other fields (such as philosophy), philosophical arguments in mathematics must be separated from ordinary mathematical deduction. Thus, the demonstration has to be formally verifiable relative to the accepted or specified axioms and the system is reduced to a recursive one. Every recursive system can be rephrased to be a recursive system of axioms with first order logic, so without loss of generality, we can limit ourselves to first order logic and recursive axiom schemes.

To increase expressive power, some had proposed the use of infinite formulas. However, human direct access to infinite formulas is unlikely in foreseeable future, so one has to work with finite representation of infinite formulas. To avoid ambiguity, the rules for representation must be formally defined, so the language is reduced to a formal language with finite formulas. Even if humans will gain an ability to work with infinite formulas, one might just as well work with sets that code the formulas and the satisfaction relation, essentially reducing the system to finite formulas plus direct access to some sets.

The expressive power of a language is determined by the objects that are described by the language and the properties expressible in the language. It turns out that in mathematics as we know it, all mathematical objects are sets or are reducible to sets, more precisely (if there are other kinds of sets) extensional well-founded sets built from the empty set. A discovery of a mathematical object that is not reducible to sets would be revolutionary and unlikely. The (first order) language of set theory consists of first order logic and the binary relation of membership (equality is definable using extensionality). That language is the de facto formal language of mathematics, and (in cases of doubt) meaningfulness of other formal languages is usually established by reduction to that language. Thus, the problem of inexpressibility in mathematics is essentially the problem of how to increase expressiveness of the language of set theory, and the solution is finding and formalizing properties of sets that are not expressible in that language. The reasons above also show that incompleteness in mathematics is incompleteness in set theory.

Weak extensions to the language can be formally defined by a statement in the extended language such that the statement uniquely fixes the extension. Stronger extensions are not defined formally but are given in the human language. One can specify the intent, argue and hope that exactly one extension satisfies the intent, and use informal reasoning to derive results. For acceptance, one needs to axiomatize enough of the intent and the results so that basic theory can be developed formally.

A large class of true statements, namely reflection principles, forms a natural hierarchy. A reflection principle can usually be represented in the form that such and such natural theory has a sufficiently well behaved model, and statements of that form tend to be reflection principles. Reflection principles are ordered by strength and expression level. A reflection principle implies all reflection principles of lower strength and expression level, including the restriction of more expressive but weaker principles to lower expression level (except at the expressive level of large cardinals when it is unstated whether the cumulative hierarchy reaches the cardinals). A sufficiently strong reflection principle resolves all basic incompleteness at its expression level (at least if the expression level is not too high; at $\Sigma^2_1$ level one would need to combine a reflection principle with a forcing statement such as the Continuum Hypothesis). For example, projective determinacy resolves the basic incompleteness for $\Sigma^1_{<\omega}$ formulas. The incompleteness addressed in this paper is at high expression level rather than high strength and is addressed using reflection principles.

# Almost Self-Referential Formulas

## The First Extension

Let Tr be the satisfaction relation for one parameter formulas of ordinary set theory. By Tarski's theorem of undefinability of truth, Tr is not definable in ordinary set theory, so let us consider adding it to the language. Almost self-referential logic can be seen as a generalization and iteration of this extension. We have to show that Tr is well-defined. Necessarily, the demonstration is a philosophical, not a formal mathematical one. However, the philosophical argument corresponds to a mathematical proof that the extension is well-defined when the domain is a set. If Tr is not well-defined, then there is a code n such that Tr(n) is not well-defined, but Tr(n)↔φ for a formula φ in ordinary set theory and so must be well-defined since φ is well-defined; φ is well-defined since 'there exist', 'or', 'not', 'in', and 'set' are (assumed to be) well-defined.

To axiomatize this extension, we add axioms to the effect that Tr is a satisfaction relation and the replacement axiom schema for formulas involving Tr. The replacement schema is true for the same reasons that ordinary replacement is true. The extension of ZFC is weaker than ZFC + there is an inaccessible cardinal. However, the extension resolves some incompleteness at a very high expression level, so some of its consequences (in the ordinary language) are not provable from ZFC + there is a proper class of supercompact cardinals if the later is consistent.

## Almost Self-Referential Formulas

Simple self-reference is not allowed since otherwise the utterance "This statement is false" would lead to a paradox. On the other hand, properties can be defined inductively (for example, $2^n$ is $2^{n-1}$ if n>0 and 1 if n=0). Such reference is possible because the parameters are varied on a well-founded path. In general, variation on a non-well-founded path is forbidden since the definition $P(n) \leftrightarrow \forall m>n \ \neg P(m)$ is self-contradictory. A formalization of the allowed self-reference is:

**Definition:** Almost Self-Referential Formula
<u>Syntax:</u> Consists of the formula proper φ(x,...) and a binary relation R. The relation is an ordinary formula with markings to indicate the arguments. The formula proper is like an ordinary formula except that it may contain a special predicate symbol *P* (same arity as φ).

<u>Semantics:</u> The truth-value of φ is computed as if it is an ordinary formula, except for *P*. Let $R_w$ be the well-founded part of R: $R_w$(s, t) iff R(s, t) and there is no infinite sequence t, $t_1$, $t_2$,..., such that R($t_1$, t) and R($t_2$, $t_1$) and ... . If *P* is inside φ(t,...), then *P*(s,$a_1$,..,$a_n$) is equivalent to $R_w$(s, t) ∧ φ(s,$a_1$,..,$a_n$).

A *nested almost self-referential formula* is as above but allowing R to use a nested almost self-referential formula (but R may not refer to φ); a formula that uses a nested almost self-referential formula as a subformula is also a nested almost self-referential formula.

For example, x is a power of 2 can be denoted by (x=1 ∨ (x is even ∧ *P*(x/2)), R is '<' restricted to natural numbers).

Consider the language of set theory extended with the satisfaction relation for almost self-referential formulas. To see that the extension is meaningful, suppose that the relation (and hence a particular almost self-referential formula) is not well-defined. Since $R_w$ is well-defined and well-founded, there must be a minimal x for which the formula is not well-defined. Since the formula proper is like an ordinary formula, that can only happen if an invocation of *P*, say with parameter y, is not well-defined. If not $R_w$(y, x), then *P*(y) is false and hence well-defined. Otherwise, $R_w$(y, x) holds, so by minimality of x, the almost self-referential formula with y in place of x, and hence *P*(y), is well-defined. Therefore, the extension is meaningful. (Note that if the domain (including quantification domain) of φ were a set, the argument (properly formalized) would be provable, but for V it is meant to be an intuitive argument. Whether something is well-defined is not in general well-defined, so formally existence of "a minimal x for which the formula is not well-defined" (above) uses an ontological assumption.) The expressive power will be analyzed, and an axiomatization presented in the next section (subsection Axiomatization of Extensions).

The extension with almost self-referential formulas is also meaningful in other contexts provided the base theory includes basic set theory and the notion of being well-founded. For example, it is meaningful (and expressible in set theory) for second order arithmetic.

A different extension is through infinitary logic, specifically $L_{\infty,\omega}$. An infinitary formula is like an ordinary formula except that disjunctions (and conjunctions) may consist of infinitely many terms (such as x=1 or x=2 or x=3 or ...). A good background on infinitary logic can be found in [3]. The extension adds satisfaction relation for one-parameter infinitary formulas. Satisfaction relation for infinitary sentences is just as expressive but is less convenient for analyzing restrictions of infinitary logic.

## Extensions

Expressive power of almost self-referential formulas can be increased by increasing the rank of the well-founded relation by approaching self-reference. This can be done by adding a second well-founded relation and variable. In that extension, self-reference requires decreasing (second variable, first variable) in lexicographical order. Reference to the original formula in the first relation is now allowed provided that the second variable decreases. Nested almost self-referential formulas correspond to the second relation having finite rank; nested almost self-referential formulas without nesting in R correspond to ordinary almost self-referential formulas.

The system can be further extended by increasing the number of relations. In fact, the number of the relations can be made unlimited, with the relations indexed by a special well-founded relation (that may not refer to φ). The relations would be coded on paper by a tertiary relation, the first argument of which is the relation index. The self-referential calls use finitely many parameters, each parameter with the index of the corresponding relation; the indices must be in the decreasing order. Each call (or level of recursion) decreases the parameter sequence in lexicographical order; and calls placed from the $j^{\text{th}}$ relation must involve a decrease in index(es) above $j$. The mechanics of these extensions resemble ordinal notation systems, and this particular extension resembles the large Veblen ordinal.

Further extensions can of course be made. However, without an overarching principle, they would only add complexity and philosophical doubt while providing less and less utility.

## Expressive Power of the Extensions

### Constructible Hierarchy Above V

The expressive power of almost self-referential logic is best explored using transitive models. Consider a transitive model V — meant to represent the universe — of ZFC or ZFC minus power set. We assume that V satisfies the replacement schema for formulas in the extension we are considering (this is not needed for the theorem). $L_α(V)$ is level α of the constructible hierarchy built above V. The expressive power of a logic is measured by what subsets of V are definable in V using that logic, allowing parameters. For infinitary logic, we only include formulas that are in V. We now state the theorem.

**Theorem:** Let V be a transitive model of ZFC (or ZFC\P). For every non-zero ordinal α in V and set S⊂V, the following are equivalent:

1. $S \in L_α(V)$
2. S is expressible in V using infinitary logic with formulas restricted to rank less than ωα (the formula is in V, and may use parameters in V).
3. S is expressible in almost self-referential logic in V (allowing parameters in V) with the well-founded relation restricted to a rank less than ωα. (One formulation is S={s:φ(a,s)} where a∈V and φ is an almost self-referential formula (with the self-reference acting on *a*).)

The proof is a routine induction on α.

For the well-founded relations used in almost self-referential formulas, the most important property is rank. For natural extensions and fragments, expressive power is linearly ordered by the allowed rank. Ranks above Ord (the height of V) lead to expressiveness beyond infinitary logic. (The length or rank of proper class relations can be compared inductively using almost self-referential logic.)

The satisfaction relation of almost self-referential logic (and its above-mentioned multivariable extensions) is definable in second order logic about V. The definition can be chosen to be absolute for every transitive extension of V (that includes V as a set) satisfying enough set theory, such as ZF\P with separation and replacement limited to $Σ_1$ formulas. In fact, KP appears to suffice, with the extended replacement in V giving sufficient absoluteness for appropriate S⊂V being well-founded. Otherwise, KPi (KP plus every set is inside an admissible set) would suffice.

## Axiomatization of Extensions

Observe that of the extensions considered, if extension A is more expressive than extension B, and if V satisfies replacement schema for A (and power set), then V also satisfies, "there is a level of the cumulative hierarchy, $V_κ$, that satisfies all true statements expressible in extension B." Therefore, without fear of falsehood, we can axiomatize an extension A by augmenting ZFC as follows:

*Axiomatization Template* (using ZFC as base):
*Variant 1:* (schema over φ and ψ in the extended language) If every $V_κ$ satisfying replacement axiom instance ψ satisfies φ, then φ.
There are natural variants:
*Variant 2:* (schema) If ZFC proves that "every $V_κ$ satisfying replacement for all formulas in A also satisfies φ", then φ.
*Variant 3:* (schema) If every $V_κ$ satisfying replacement for all formulas in A also satisfies φ, then φ.
Large cardinal variants (C is a large cardinal axiom; C(κ) implies that κ is inaccessible):
*Variant 4:* (schema) If ZFC proves "every $V_κ$ with C(κ) satisfies φ", then φ.
*Variant 5:* (schema) If every $V_κ$ with C(κ) satisfies φ, then φ.

Such an axiomatization resolves the need for deductive apparatus for A. For typical A, variant 3 is slightly stronger than variant 2, which in turn is slightly stronger than variant 1.

Variants 3 and 4 allows formalization of the idea that the class of ordinals, Ord, behaves like a large cardinal. Such formalization extends the strength of the large cardinal to the highest available level of expression. Compared to variant 4, variant 5 increases the strength but also require greater ontological commitment; variant 4 does not require the large cardinal to exist in V (and a slight weakening of variant 4 is equiconsistent with ZFC + ∃κ C(κ)). However, variant 4 with a mild strengthening of the large cardinal C implies variant 5 with unstrengthened C; for the extensions above (using almost self-referential formulas and their extensions), it suffices to replace large cardinal with regular limit of stationary many of these large cardinals.

There is a close connection between the axiomatization for almost self-referential formulas, and NBG with determinacy for clopen class-length games. The later is described in [4], where the determinacy is shown equivalent to ETR: Elementary transfinite recursion over well-founded class relations. ETR essentially corresponds to almost self-referential formulas. To get what appears to be the exact strength (and the same theorems for statements expressible in the common language), add finite iterations of the satisfaction relation for almost self-referential formulas, or work with the formulas directly but allow nested almost self-referential formulas; use variant 1 of the axiomatization above.

## A General Extension

The general idea behind the almost self-referential extensions is that restricted self-reference is meaningful when there is a guarantee that there are no contradictions and ambiguities. That idea can be formalized in a general way. Fix a Gödel numbering $\{\varphi_0, \varphi_1, \varphi_2 ...\}$ of sentences of set theory with a predicate symbol S. Call a set M good if M is $V_\kappa$ with κ having uncountable cofinality, and call n good if ∀(good M) ∃!S⊂M (M,∈,S)⊨φ. The expressiveness of the extension can be varied by changing the conditions on "good". Uncountable cofinality was used to ensure that well-foundness for S is expressible; it is unclear if that is important. Let a predicate $P_n$ be such that $P_n(x) \Leftrightarrow$ (n is good) ∧ $\varphi_n[S/P_n]$. Define P so that ∀n∀x (P(n, x)↔$P_n$(x)). By reflection, the uniqueness of S for all sufficiently well-behaved models implies uniqueness of $P_n$, so P is unique in V. This extension is more expressive than almost self-referential logic and its extensions, but any reasonable variation is less expressive than $\Sigma^1_1$ formulas of second order logic about V. The later claim can be proved in (for example) KM (NBG with full impredicative class comprehension), assuming that the condition on good is sound (in that P is unique) and expressible in (V,∈); KM more than suffices to prove soundness of the above extension.

A less expressive extension is the least fixed point logic. That logic includes the least fixed point operator, lfp(φ) is the least (i.e. the least under '⊂') X such that X = φ(X), where X only occurs positively in φ (φ may have additional free variables). The logic has no other second order quantification, so X and lfp can only be used as in X(x) and lfp(φ)(x). Even a single lfp is more expressive than almost self-referential logic (and its above-mentioned multivariable extension), at least for models closed under countable sequences and satisfying basic set theory. Some variations of fixed point logic are described in [5].

## A Hierarchy of Large Cardinals

The replacement schema can be viewed as a statement that the class of ordinals, Ord, is inaccessible. We can formalize this as follows: $L_\alpha(V)$ satisfies "Ord(V) is inaccessible" iff the corresponding logic about V satisfies the replacement schema. The correspondence (between α and the choice of logic) is that of expressive power and includes the correspondence stated in the theorem of the previous section. A similar correspondence also extends to Mahlo cardinals by replacing inaccessible with Mahlo and requiring that V satisfies the axiom schema MAH: Every continuous class function on ordinals has a regular fix point (each instance of the schema is for a separate formula coding the function; all formulas in the logic considered can be used). Equivalently, MAH states that for every one-parameter formula φ, there is an inaccessible ordinal κ such that ∀x∈$V_\kappa$ (φ(x)↔"φ(x) holds in $V_\kappa$").

The correspondence of models can be used to yield the following equiconsistency result.

**Theorem:**

1. ZFC minus infinity is equiconsistent with rudimentary set theory plus there is an infinite ordinal.
2. ZFC minus power set is equiconsistent with rudimentary set theory plus there is an uncountable ordinal.
3. ZFC is equiconsistent with rudimentary set theory plus there is an inaccessible cardinal.
4. ZFC + MAH is equiconsistent with rudimentary set theory plus there exists a Mahlo cardinal.

> Note: More generally, ZFC + (schema) {there is φ-correct large cardinal : φ is a formula} is equiconsistent with rudimentary set theory + there is a regular limit of stationary many of these cardinals.

Rudimentary set theory is meant to be the weakest set theory in which some basic things can be done. Levels of Jensen hierarchy for L satisfy it. It is not clear how strong the theory should be, but for the theorem the following version/axiomatization works: extensionality, foundation, empty set, pairing, union, existence of transitive closure, existence of the set of all sets with transitive closure less numerous than a given set, and bounded quantifier separation. The theorem generalizes to more expressive logics by replacing rudimentary set theory with an appropriate mild strengthening: ZF minus power set corresponds to second order logic.

**Measurability and Ord:** We now argue that Ord is like a measurable cardinal in the sense that only true statements can be proven in ZFC to hold in every $V_κ$ with measurable κ. Take a mental image of the universe, and try extending it with new cardinals. By definition, the universe has all sets it may possibly contain, so there must be some reason why the new cardinals are not realized in the universe. A good reason is non-rigidity, for imagination is non-rigid while mathematical reality is fixed. The required non-rigidity seems to be (at least) the existence of a nontrivial elementary embedding of V into M if the cardinals existed. Such embedding raises the question as to why the universe is V rather than M; metaphorically speaking, such embedding implies that V is not held together by its structure. The hierarchy of cardinals below measurable is canonical and well understood, and a number of arguments exists for various such cardinals. Thus, Ord is measurable in the extension of the mental image, and since such images (at least in theory) can be made with arbitrary accuracy, Ord can be made measurable in the above-mentioned sense.

# Proper Classes

A proper class would be a collection of objects that is too large to be a set. The problem is that a set is any collection of objects. The cumulative hierarchy is without end, so if proper classes did exist, we could form a layer on top of them (and obtain a hyper-class of all classes of ordinals), and after adding sufficiently many layers, discover that the previous notion of set was too restrictive and that what was considered to be proper classes are actually sets. However, the word "proper class" is frequently used in set theory and category theory, so some treatment is needed. There are three alternatives.

* One can treat proper classes as syntactical objects and sets as sets and work in a theory conservative over set theory. For example, NBG (which includes global choice) is conservative over ZFC. However, one can do much better using the strong reflection principles proposed: Ord would be a large cardinal; the classes would satisfy full comprehension, and one can even build hyper-classes and higher types with appropriate comprehension.

* One can consider only sets whose rank is smaller than a certain ordinal κ, and treat collections of these sets as classes. The ordinal should be sufficiently large to include all objects under consideration, and $V_κ$ should satisfy enough correctness so that relevant statements, if true in $V_κ$, are true in V. Strong reflection principles allow one to assume that $V_κ$ is sufficiently closed (at the least, κ should be inaccessible, but

weak compactness, total indescribability, and more are desirable). Almost self-referential logic (with the corresponding enhancement of replacement) allows selection of κ with unprecedented degree of correctness: $V_κ$ could be an elementary substructure of V and some second order statements about $V_κ$ transfer (in a truth preserving way) into almost self-referential statements about V.

\* One can treat classes semantically as descriptions, and for a set x and proper class X, treat x ∈ X as "set x satisfies description X". Almost self-referential logic allows much richer descriptions, and hence much stronger closure properties. For example, if the descriptions are formulas in infinitary logic, then a union/disjunction of (set size many) proper classes is a proper class.

Note: If one believes in proper classes, then one can use almost self-referential logic to make (what one would believe are) more expressive statements about proper classes.

There may be a fourth alternative. . .

# Higher Order Set Theory

## Unsuccessful Attempts

Simple attempts at higher order set theory appear trivial or meaningless. Allowing infinite blocks of quantifiers in infinitary formulas does not increase expressive power of set theory because in place of "∃$x_1$ ∃$x_2$ ..., such that ..." one can use "there is X=($x_1$, $x_2$, $x_3$ ...) such that ...", and use the $n^{th}$ element of X in place of $x_n$.

Infinite alterations of quantifiers are not generally meaningful because some games are not determined. However, if one limits the universe to real numbers, then (under sufficient large cardinal assumptions) ∃$x_1$ ∀$x_2$ ∃$x_3$ ∀$x_4$ ... f($x_1$, $x_2$, $x_3$, $x_4$ ...) is meaningful and inclusion of such formulas increases the expressive power to that of L(**R**) and beyond.

Second order logic about sets does not increase expressive power either (unless one works in a system where the power set axiom is false), as second order statements about S are first order statements about the power set of S. Second order logic about V appears meaningless because (in the required sense) there are no proper classes. However, there may be an alternative.

## Reflective Ordinals

Imagine reaching ordinals of higher and higher large cardinal and reflecting properties. One reaches inaccessible, Mahlo, weakly compact, and other cardinals. One reaches $Σ_2$ correct ordinals, $Σ_3$ correct ordinals, levels of the cumulative hierarchy that are elementary substructures of V, and so on. A natural hypothesis is that as one reaches higher and higher levels of reflection, the properties of the ordinals converge. Specifically,

**Convergence Hypothesis:** If α and β are ordinals with sufficiently strong reflecting properties, then φ(S, α) ↔ φ(S, β) whenever rank(S) < min(α, β) and φ is a first order formula of set theory with two free variables.

The hypothesis is a hypothesis in a philosophical rather than a formal mathematical sense.

If the convergence hypothesis holds, let us call an ordinal α *reflective*, and denote it by R(α) when for every S∈$V_α$ and for every first order formula in set theory and with two free variables φ, φ(S, α) iff φ(S, β) for every β>rank(S) with sufficiently strong reflecting properties. If zero sharp exists, then the convergence hypothesis holds in L, with every Silver indiscernible having sufficiently strong reflecting properties. For every indiscernible α, $R^L$ restricted to α is constructible and is uniformly definable in L from (α, β) where β is any indiscernible above α.

To get a reasonable axiomatization of R, recall the imagination of Ord being measurable. For an extender E on κ, let j: V→M be the corresponding elementary embedding. κ represents Ord and hence can consistently act as a reflective ordinal in j($V_κ$). Define $R_E$ by $R_E$(α) iff α∈κ and j($V_κ$) satisfies φ(S, α)↔φ(S, κ) for every S∈$V_α$ and every formula φ (in the first order language of set theory) with two free variables. $R_E$ is unbounded in κ and has measure one, and every member of $R_E$ is totally indescribable, totally ineffable, Ramsey, and so on.

Let us augment the language of set theory with a unary predicate R, and add to ZFC the following axiom schema: If ZFC proves that "for every measurable cardinal κ and every extender E with critical point κ, sentence φ is true in ($V_κ$, ∈, $R_E$)", then φ.

The axiomatization is consistent relative to an inaccessible cardinal above a measurable one. The axiomatization can be weakened to the level of indescribable cardinals by using embeddings of only a subset of V. [Update (Dec 16, 2016): One needs a strong version of indescribability; see "Reflective Cardinals" for details.] The axiomatization is reasonably complete (at least with respect to large cardinal structure) if we assume that there are no measurable cardinals. However, to be non-restrictive, we make no such assumption, and we allow the embedding to be witnessed by an extender rather than only by an ultrafilter. [Update/correction (Dec 16, 2016): An ultrafilter would suffice since E and the ultrafilter on P(κ) corresponding to E would lead to the same R.] If there are measurable cardinals in V, then an appropriate strengthening of the axiomatic system above may be to require κ be measurable in M, and again, if we assume that there are no (measurable) cardinals of order 2 in V, we get a reasonably complete theory. Axiomatization of reflective ordinals should be guided by the fact that reflective ordinals satisfy all true large cardinal properties that are realized in V and are expressible in first order set theory.

By the reflection principle, in so far as higher order statements in set theory are meaningful, they are true about V iff they are true about $V_κ$ for reflective ordinals κ. For example, an infinitary statement S is true iff for all α>rank(S), R(α) implies S is true in $V_α$. However, we can turn this around and use reflective ordinals to *define* semantics for higher order set theory. If φ is a formula in higher order set theory, then φ($s_1$, $s_2$, ..., $s_n$) (where each $s_i$ is a set) is true iff for every reflective ordinal α > rank({$s_1$, $s_2$, ..., $s_n$}), $φ_1$(α, $s_1$, $s_2$, ..., $s_n$) holds, where $φ_1$ is the translation of φ accomplished in the same way as translation of higher order formulas about $V_α$ into formulas of set theory. The metaphysical semantics of proper classes depends on the reflective ordinal chosen; particular classes can be viewed as descriptions; and particular predicates on classes—as template descriptions, which are materialized once an appropriate notion of proper classes is chosen.

Thus, if R is well-defined, it supersedes almost self-referential formulas. Proper classes are finally given proper semantic treatment. The expressiveness of the new language is so great that it eliminates the need for using almost self-referential formulas. A skeptic may well claim that the notion of a reflective ordinal is mysticism, but even if it is, the formalization above is a valid and interesting formal system, and, in the author's opinion, is the most natural way to formalize higher order set theory.

## Extensions

The language of set theory with reflective ordinals satisfies current mathematical needs, and it is too early to accept further extensions. However, the following two related questions may make further study important:

- Is there a definable set that is not ordinal definable?
- Is it possible to define (the real number coding) first order set theory locally?

Locally means without invoking large cardinals; definitions in $V_{\omega+\omega}$ are local, but definitions invoking the supremum of ordinals definable in first order set theory are not. An affirmative answer to the second question could radically revise our understanding of the mathematical universe. It would imply existence of properties inherently more expressive than those currently accepted.

The notion of a reflective ordinal has a light-face weakening: An ordinal *s* would be called light-face reflective if the theory of (V, ∈, *s*) would remain the same had *s* been replaced by a reflective ordinal. An extension to set theory could be made simply by augmenting symbol *s* for the least light-face reflective ordinal, but the system in the paragraphs above is preferred because many more predicates (such as the satisfaction relation for first order set theory) are expressible in that system.

Recall that a reflective ordinal κ is an ordinal sufficiently large relative to smaller ordinals for the theory (V, ∈, κ) with parameters in $V_\kappa$ to be "correct". The idea of reflective ordinals can be iterated and applied to other extensions of set theory. An ordinal κ is reflective in degree α+1 iff the theory (V, ∈, κ, R) with parameters in $V_\kappa$ is correct, where R is the predicate for reflectiveness in degree α. For a limit ordinal α, reflectiveness of degree α can be defined as reflectiveness in all degrees less than α; every ordinal is reflective in degree 0. In L, the indiscernibles are reflective in every finite degree, but non-existence of $0^\#$ in L prevents formalization except for fixed finite degrees.

Higher degree reflective ordinals can be axiomatized inductively using elementary embeddings analogously to the axiomatization of reflective ordinals. In the models, the critical point κ (representing Ord) can be reflective in degree α+1 in $j(V_\kappa)$; j can extend degree α reflectiveness in $V_\kappa$ until j(κ); and α+1 degree reflectiveness of κ can be used to define such reflectiveness below κ. (However, it is unclear whether using a single elementary embedding for all degrees of reflectiveness is correct. [Update/correction (Dec 16, 2016): Using a single embedding is correct.])

A different extension is through reflective sequences. The predicate for reflective sequences (call it P) is defined by (if what is below counts as a definition):

1. Every set of ordinals that excludes its limit points and has sufficiently strong reflecting properties (for its order type) satisfies P.
2. P(S) iff S is a set of ordinals and for every set of ordinals T of the same order

type satisfying P,
$\forall x \in V_{\min(S \cup T)}$ Theory(V, ∈, S, x) = Theory(V, ∈, T, x).
[Update (Dec 16, 2016): There are variations on the definition in how to handle x between min(S) and sup(S).]

Clearly, the order type of a reflective sequence is less than its minimum. However, there may be a stronger extension defined through correctness of Theory(V, ∈, S) with parameters in $V_{\min(S)}$ , where the set of ordinals S excludes its limit points, has order type with sufficient reflecting properties and larger than any element of S, and with every member of S having sufficient reflecting properties for the sequence. Conceivably, the stronger extension has an elegant weakening (which is more expressive than P), namely the predicate S that acts as the maximal proper class sequence of ordinals with sufficient reflecting properties for Theory(V, ∈, S\rank({x}), x) to be correct for every set x.

The predicate for reflective ordinals of a finite degree n>1 has slightly greater expressive power than the predicate for reflective n-tuples of ordinals (this also applies to L). For example, every {reflective ordinal of second degree, a larger reflective ordinal} is a reflective pair. However, infinite reflective sequences have greater expressive power than iterations of reflectiveness: If S is such a sequence, then for every ordinal α, reflectiveness of degree α restricted to min(S) is definable from S and α. It is unclear if infinite reflective sequences exist and how to axiomatize them.

For a fixed natural number n and ordinal α, reflectiveness for n-tuples below α is ordinal definable. However, could it be that a real number is ordinal definable iff it is definable from a reflective n-tuple of ordinals?

The way to proceed may be study analogues of these notions for hereditarily countable sets and other well-understood theories. Perhaps, such studies will lead to a local definition of first order set theory.